%% file: preprint.tex
\documentclass[%
  letterpaper,
  twocolumn,
  colorlinks,
]{preprint}

\usepackage[noadjust]{cite}

\usepackage[english]{babel}

\usepackage{amsmath}
\usepackage{amsfonts}
\usepackage{amssymb}
\usepackage{amsthm}
\usepackage{amscd}
\usepackage{bbm}
\usepackage{dsfont}

\usepackage{graphicx}

\usepackage{tikz}
\usepackage{pgfplots}
  \pgfplotsset{compat = 1.17}
  \pgfkeys{/pgf/number format/.cd, 1000 sep = {\,}}
  \usepgfplotslibrary{external}
  \tikzset{external/system call = {%
    pdflatex \tikzexternalcheckshellescape
      -halt-on-error
      -interaction=batchmode
      -jobname "\image" "\texsource"}}
  \tikzexternaldisable

\usepackage{subcaption}

\newcommand{\trans}{\ensuremath{\mkern-1.5mu\mathsf{T}}}
\newcommand{\herm}{\ensuremath{\mathsf{H}}}
\DeclareMathOperator{\mspan}{span}
\DeclareMathOperator{\diag}{diag}

\newcommand{\C}{\ensuremath{\mathbb{C}}}
\newcommand{\Cn}{\ensuremath{\C^{n}}}
\newcommand{\Cp}{\ensuremath{\C^{p}}}
\newcommand{\Cnn}{\ensuremath{\C^{n \times n}}}
\newcommand{\Cpn}{\ensuremath{\C^{p \times n}}}
\newcommand{\Cr}{\ensuremath{\C^{r}}}
\newcommand{\Cnr}{\ensuremath{\C^{n \times r}}}
\newcommand{\Crr}{\ensuremath{\C^{r \times r}}}

\newcommand{\R}{\ensuremath{\mathbb{R}}}
\newcommand{\Rn}{\ensuremath{\R^{n}}}
\newcommand{\Rnn}{\ensuremath{\R^{n \times n}}}

\renewcommand{\i}{\ensuremath{\mathfrak{i}}}
\newcommand{\Sys}{\ensuremath{\Sigma}}
\newcommand{\LinSys}{\ensuremath{\Sys_{\operatorname{LO}}}}
\newcommand{\SoSys}{\ensuremath{\Sys_{\operatorname{SO,LO}}}}
\newcommand{\QOSys}{\ensuremath{\Sys_{\operatorname{QO}}}}
\newcommand{\QOSystd}{\ensuremath{\Sys_{\operatorname{QO}}^{\operatorname{t}}}}
\newcommand{\QOSysred}{\ensuremath{\widehat{\Sys}_{\operatorname{QO}}}}
\newcommand{\relerr}{\ensuremath{\operatorname{relerr}}}
\newcommand{\nso}{\ensuremath{n_{\operatorname{so}}}}

\newcommand{\BM}{\ensuremath{\boldsymbol{M}}}
\newcommand{\BD}{\ensuremath{\boldsymbol{D}}}
\newcommand{\BK}{\ensuremath{\boldsymbol{K}}}
\newcommand{\BA}{\ensuremath{\boldsymbol{A}}}
\newcommand{\BC}{\ensuremath{\boldsymbol{C}}}
\newcommand{\BE}{\ensuremath{\boldsymbol{E}}}
\newcommand{\BG}{\ensuremath{\boldsymbol{G}}}
\newcommand{\BQ}{\ensuremath{\boldsymbol{Q}}}
\newcommand{\BI}{\ensuremath{\boldsymbol{I}}}
\newcommand{\BV}{\ensuremath{\boldsymbol{V}}}
\newcommand{\BW}{\ensuremath{\boldsymbol{W}}}
\newcommand{\Bb}{\ensuremath{\boldsymbol{b}}}
\newcommand{\Bc}{\ensuremath{\boldsymbol{c}}}
\newcommand{\Be}{\ensuremath{\boldsymbol{e}}}
\newcommand{\Bg}{\ensuremath{\boldsymbol{g}}}
\newcommand{\Bp}{\ensuremath{\boldsymbol{p}}}
\newcommand{\Bv}{\ensuremath{\boldsymbol{v}}}
\newcommand{\Bw}{\ensuremath{\boldsymbol{w}}}
\newcommand{\Bx}{\ensuremath{\boldsymbol{x}}}
\newcommand{\Bz}{\ensuremath{\boldsymbol{z}}}
\newcommand{\Bzero}{\ensuremath{\boldsymbol{0}}}

\newcommand{\BVr}{\ensuremath{\BV}} 
\newcommand{\BWr}{\ensuremath{\BW}}

\newcommand{\BAr}{\ensuremath{\skew5\widehat{\BA}}}
\newcommand{\BCr}{\ensuremath{\skew2\widehat{\BC}}}
\newcommand{\BEr}{\ensuremath{\skew4\widehat{\BE}}}
\newcommand{\BGr}{\ensuremath{\skew2\widehat{\BG}}}
\newcommand{\BQr}{\ensuremath{\skew2\widehat{\BQ}}}
\newcommand{\Bxr}{\ensuremath{\skew1\widehat{\Bx}}}
\newcommand{\Bbr}{\ensuremath{\skew0\widehat{\Bb}}}
\newcommand{\Hr}{\ensuremath{\skew2\widehat{H}}}
\newcommand{\yr}{\ensuremath{\widehat{y}}}
\newcommand{\qr}{\ensuremath{\skew{-.5}\widehat{q}}}

\newcommand{\Bxs}{\ensuremath{\check{\Bx}}}
\newcommand{\xs}{\ensuremath{\check{x}}}

\newcommand{\Bxt}{\ensuremath{\Bx_{\operatorname{t}}}}
\newcommand{\Bdxt}{\ensuremath{\dot{\Bx}_{\operatorname{t}}}}
\newcommand{\Ht}{\ensuremath{H_{\operatorname{t}}}}
\newcommand{\Hrt}{\ensuremath{\skew2\widehat{H}_{\operatorname{t}}}}
\newcommand{\ut}{\ensuremath{u_{\operatorname{t}}}}
\newcommand{\yt}{\ensuremath{y_{\operatorname{t}}}}

\newcommand{\CH}{\ensuremath{\mathcal{H}}}
\newcommand{\CK}{\ensuremath{\mathcal{K}}}
\newcommand{\CL}{\ensuremath{\mathcal{L}}}

\newcommand{\CKr}{\ensuremath{\skew2\widehat{\mathcal{K}}}}

\newcommand{\irka}{\ensuremath{\mathsf{LQO\mbox{-}IRKA}}}
\newcommand{\intinfv}{\ensuremath{\mathsf{Int_{\infty, V}}}}
\newcommand{\intinfvw}{\ensuremath{\mathsf{Int_{\infty, VW}}}}
\newcommand{\intavgv}{\ensuremath{\mathsf{Int_{avg, V}}}}
\newcommand{\intavgvw}{\ensuremath{\mathsf{Int_{avg, VW}}}}

\theoremstyle{plain}\newtheorem{theorem}{Theorem}

\definecolor{matlabblue}{HTML}{0072BD}
\definecolor{matlaborange}{HTML}{D95319}
\definecolor{matlabyellow}{HTML}{EDB120}
\definecolor{matlabpurple}{HTML}{7E2F8E}
\definecolor{matlabgreen}{HTML}{77AC30}
\definecolor{matlablightblue}{HTML}{4DBEEE}
\definecolor{matlabred}{HTML}{A2142F}

\tikzstyle{sline} = [
  solid,
  line width = 1.5pt
]

\pgfplotscreateplotcyclelist{plotlist}{
  {black, sline},
  {matlabblue, sline},
  {matlaborange, sline, dashed},
  {matlabpurple, sline, dashed},
  {matlabgreen, sline, dotted},
  {matlabred, sline, dotted}
}

\newcommand{\plotfontsize}{\footnotesize}

\begin{document}
  

\title{Interpolatory model reduction of dynamical systems with root mean
  squared error}
  
\author[$\ast$]{Sean Reiter}
\affil[$\ast$]{Department of Mathematics, Virginia Tech,
  Blacksburg, VA 24061, USA.\authorcr
  \email{seanr7@vt.edu}, \orcid{0000-0002-7510-1530}}
  
\author[$\dagger$]{Steffen W. R. Werner}
\affil[$\dagger$]{Department of Mathematics and
  Division of Computational Modeling and Data Analytics,
  Academy of Data Science, Virginia Tech,
  Blacksburg, VA 24061, USA.\authorcr
  \email{steffen.werner@vt.edu}, \orcid{0000-0003-1667-4862}}
  
\shorttitle{Interpolatory model order reduction for root mean squared error
  systems}
\shortauthor{S. Reiter, S. W. R. Werner}
\shortdate{2025-04-18}
\shortinstitute{}
  

  
\abstract{%
  The root mean squared error is an important measure used in a variety of
  applications such as structural dynamics and acoustics to model averaged
  deviations from standard behavior.
  For large-scale systems, simulations of this quantity quickly become
  computationally prohibitive. 
  Classical model order reduction techniques attempt to resolve this issue via
  the construction of surrogate models that emulate the root mean squared error
  measure using an intermediate linear system.
  However, this approach requires a potentially large number of linear outputs,
  which can be disadvantageous in the design of reduced-order models.
  In this work, we consider directly the root mean squared error as the
  quantity of interest using the concept of quadratic-output models and
  propose several new model reduction techniques for the construction of
  appropriate surrogates.
  We test the proposed methods on a model for the vibrational response of a
  plate with tuned vibration absorbers.
}

\novelty{}

\maketitle

\section{Introduction}%
\label{sec:intro}

The modeling and numerical simulation of large-scale dynamical systems are
powerful tools for deciphering the behaviour of complex physical phenomena.
Such large-scale systems typically arise from the demand for highly accurate
predictions.
The systems that we consider in this work are described in the frequency
(Laplace) domain by linear frequency-dependent algebraic equations of the
form
\begin{equation} \label{eqn:state}
  (s\BE - \BA) \Bx(s) = \Bb u(s),
\end{equation}
where $\BE,\BA \in \Cnn$ and $\Bb \in \Cn$.
The internal state is given by $\Bx\colon \C \to \Cn$, and
$u\colon \C \to \C$ is the external input.
For simplicity, we assume throughout this paper that $\BE$ is a nonsingular
matrix.
Systems of equations of the form~\cref{eqn:state} are used to model a
wide range of different dynamical behaviors including heat transfers, fluid
dynamics or the deformation of mechanical structures; see,
for example,~\cite{BenMS05, Wer21}.

In regards of applications, typically not the full state
of~\cref{eqn:state} is observed but only a limited number of quantities
of interest.
In this work, we consider the case that the
\emph{root mean squared (RMS) error} of the state $\Bx$ in~\cref{eqn:state}
is measured with respect to a point of reference $\Bxs \in \Cn$. 
That is, we are interested in the observable
$y\colon \C \to \C$ defined as
\begin{equation} \label{eqn:rms}
  y(s) = \sqrt{\sum_{k = 1}^{n} \lvert q_{k} \rvert^{2} \,
    \lvert x_{k}(s) - \xs_{k} \rvert^{2}},
\end{equation}
where $x_k(s)$ is the $k$-th component of the state vector $\Bx(s)$,
$\xs_{k}$ is the $k$-th component of the reference point $\Bxs$,
and the scalars $q_{k} \in \C$ are weights.
Systems that consider~\cref{eqn:rms} as the observable quantity of interest
arise in a variety of different applications, for example, in the study of
structural dynamics and vibro-acoustic
systems~\cite{VanVLetal12, AumW23}, wherein one might be interested
in the average spatial deformation or displacement of a given surface.
Other examples for the use of~\cref{eqn:rms} include observables relating
to a system's internal energy~\cite{Pul23}, the variance of a random variable
in stochastic models~\cite{PulA19}, and approximations of the cost function in
quadratic regulator problems~\cite{DiaHGetal23}.
For the simplicity of presentation, we assume without loss of generality the
reference state to be the zero state, $\Bxs = \Bzero_{n}$.
This is always possible by appropriately substituting the state $\Bx$ 
in~\cref{eqn:state}.

Classically, the RMS error~\cref{eqn:rms} of the dynamical
system~\cref{eqn:state} is simulated implicitly by introducing a linear output
operator $\BC \in \Cpn$ and working with a linear input-output system of
the form
\begin{equation} \label{eqn:linsys}
  \LinSys\colon \left\{
    \begin{aligned}
      (s\BE - \BA) \Bx(s) & = \Bb u(s), \\
      \Bz(s) & = \BC \Bx(s),
    \end{aligned}\right.
\end{equation}
with the vector-valued output signal $\Bz\colon \C \to \Cp$.
The RMS error~\cref{eqn:rms} can easily be recovered from~\cref{eqn:linsys}
by choosing
\begin{equation*}
  \BC = \begin{bmatrix} q_{i_{1}} \Be_{i_{1}} & q_{i_{2}} \Be_{i_{2}} &
    \ldots & q_{i_{p}} \Be_{i_{p}}  \end{bmatrix}^{\trans},
\end{equation*}
where $\Be_{i_{j}} \in \Cn$ is the $i_{j}$-th canonical basis vector and
$i_{1}, i_{2}, \ldots, i_{p} \subseteq \{1, \ldots, n\}$ are the indices
corresponding to the nonzero weights $q_{j}$ in~\cref{eqn:rms}.
Then, for any $s \in \C$, the RMS error $y(s)$ in~\cref{eqn:rms} is given
via the $\ell_{2}$-norm of the output of~\cref{eqn:linsys}:
\begin{equation*}
  \lVert \Bz(s) \rVert_{2}
    = \sqrt{\big( \BC \Bx(s) \big)^{\herm} \big( \BC \Bx(s) \big)}
    = y(s),
\end{equation*}
where $\Bw^{\herm}$ denotes the conjugate transpose of a vector $\Bw \in \Cn$.
The frequency domain input-to-output map of the linear
system~\cref{eqn:linsys} is fully characterized by the corresponding
\emph{transfer function} $\BG\colon  \C \to \Cp$ such that
$\Bz(s)=\BG(s) u(s)$, where
\begin{equation} \label{eqn:lintf}
  \BG(s) = \BC (s\BE - \BA)^{-1} \Bb.
\end{equation}
A drawback of using~\cref{eqn:linsys} comes in the form of a
potentially large number of outputs $p$.
This stands in contrast to the one-dimensional quantity of
interest~\cref{eqn:rms} and is typically disadvantageous in the design
of reduced-order models.

We consider here a different approach to combine~\cref{eqn:state}
with~\cref{eqn:rms}.
By introducing a single quadratic-output function, the RMS error can be
directly simulated up to the square root.
Such systems with quadratic output can be written as
\begin{equation} \label{eqn:qosys}
  \QOSys\colon \left\{
    \begin{aligned}
      (s\BE - \BA) \Bx(s) & = \Bb u(s), \\
      y(s)^{2} & = \Bx(s)^{\herm} \BQ \Bx(s),
    \end{aligned}\right.
\end{equation}
where $\BQ \in \Rnn$ is Hermitian.
Colloquially, we refer to systems of the form~\cref{eqn:qosys} as
\emph{linear quadratic-output systems}.
The RMS error of $\Bx$ is easily retrieved by defining the output matrix $\BQ$ 
in~\cref{eqn:qosys} as the diagonal matrix with the weights
$\lvert q_{k} \rvert^{2} \geq 0$  as entries, i.e.,
\begin{equation*}
  \BQ = \diag(\lvert q_{1} \rvert^{2}, \lvert q_{2} \rvert^{2}, \ldots,
    \lvert q_{n} \rvert^{2}).
\end{equation*}

In practical applications, the state dimension of~\cref{eqn:qosys} becomes
easily very large ($n \in \mathcal{O}(10^{5})$ or larger) due to the demand
for highly accurate models.
In this case, any repeated action involving the full-order
model~\cref{eqn:qosys} commands significant computational resources
such as time and memory.
A remedy to this problem is \emph{model order reduction}, which is concerned
with the construction of low-order, cheap-to-evaluate surrogate systems that
can be used in place of the full-order model.
In this work, we consider the construction of reduced-order models for
linear quadratic-output systems of the form~\cref{eqn:qosys} given as
\begin{equation} \label{eqn:qosys_red}
  \QOSysred\colon \left\{
    \begin{aligned}
      (s\BEr - \BAr) \Bxr(s) & = \Bbr u(s), \\
      \yr(s)^{2} & = \Bxr(s)^{\herm} \BQr \Bxr(s),
    \end{aligned}\right.
\end{equation}
where $\BAr,\BEr,\BQr\in\Crr$, $\Bbr\in\Cr$, $\Bxr\colon \C \to \Cr$ is the
approximate internal state and $r \ll n$.
In order to be an effective replacement, the approximate
system~\cref{eqn:qosys_red} should replicate the input-to-output behaviour
of the original large-scale system~\cref{eqn:qosys}.
In other words, given a tolerance $\tau>0$, the approximate output
$\yr\colon \C \to \C$ in~\cref{eqn:qosys_red} should satisfy
\begin{equation*}
  \lVert y - \yr \rVert \leq \tau \cdot \lVert u \rVert,
\end{equation*}
in appropriate norms and for a range of admissible inputs~$u$.

Considering the classical approach to model systems with the
output~\cref{eqn:rms} via linear multi-output systems~\cref{eqn:linsys}, one
advantage is the large variety of established model reduction procedures;
see, for example,~\cite{BenMS05, AntBG20} and the references
therein.
However, this linearization approach typically necessitates a large number of
outputs $p \approx n$.
Traditional linear system approximation techniques tend to produce poor quality
approximants and require more computational effort when $p$ is
large; see~\cite{AumW23}.
On the other hand, recently there has been a surge of interest in the model
reduction of linear quadratic-output systems formulated in the
\emph{time domain}~\cite{VanM10, VanVLetal12, PulA19, GosA19, GosG22, BenGP21,
DiaHGetal23, Pul23, ReiPGetal24, ReiPGetal23}.
These systems bear many similarities to the frequency domain linear
quadratic-output systems~\cref{eqn:qosys} that hold our interest in this work;
however, there are subtle differences between these system classes, which we
will outline later on.

In this work, we advocate for the treatment of large-scale frequency domain
problems, which model the RMS error~\cref{eqn:rms} as quantity of interest
using the linear quadratic-output system class~\cref{eqn:qosys}.
A particular example for such a system from the literature is outlined in
\Cref{sec:plate_tva}.
We will make use of some of the established theory for time domain systems,
which we review in \Cref{sec:lqo_review}, before we discuss classical
interpolation approaches for linear systems~\cref{eqn:linsys} and propose novel
extensions of this interpolation theory for the system class~\cref{eqn:qosys}
we consider in this work in \Cref{sec:interp_mor}.
To verify our theoretical results, we apply these new and extended model
reduction techniques to the practical example from \Cref{sec:plate_tva}
and report the results in \Cref{sec:numerics}.
The paper is concluded in \Cref{sec:conclusions}.

\section{Vibrations of a plate with tuned vibration absorbers}%
\label{sec:plate_tva}

\begin{figure}[t]
  \centering
  \vspace{.5\baselineskip}
  \resizebox{.9\linewidth}{!}{%
  \tikzexternalenable%
  \tikzsetnextfilename{plate}%
  \input{graphics/plate.tikz}%
  \tikzexternaldisable%
}
  \caption{Visual sketch of a plate equipped with TVAs~\cite{AumW23}.}
  \label{fig:platetva}
  \vspace{-\baselineskip}
\end{figure}
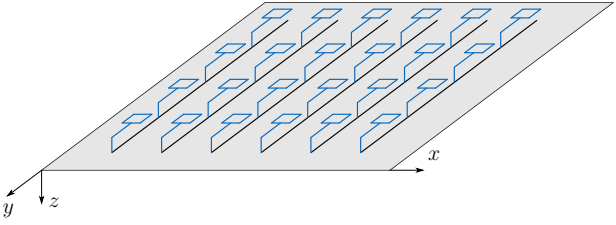

As a specific example for systems of the form~\cref{eqn:qosys}, we consider the
model of the vibrational response of a strutted plate
from~\cite[Sec.~4.2]{AumW23}.
The matrices for the computational model itself are available
at~\cite{supAumW22}.
The plate has the dimensions $0.8 \times 0.8$\,m, with a thickness of
$1$\,mm, and it consists of aluminium, with the material parameters
$E = 69$\,GPa, $\rho = 2650$\,kh\,m$^{-3}$ and $\nu = 0.22$.
The damping is assumed to be proportional (Rayleigh damping) with the scaling
parameters $\alpha = 0.01$ and $\beta = 1 \cdot 10^{-4}$.
The tuned vibration absorbers (TVAs) connected to the plate are used to reduce
the vibrational response in the frequency region about $48$\,Hz.
These TVAs are modeled as discrete mass-spring-damper systems.

Acoustical engineers are especially interested in the RMS displacement of
the plate points in $z$-direction in response to point load excitation;
cf. \Cref{fig:platetva}.
This is recovered in~\cite{AumW23} using a \emph{second-order} linear-output
system of the form
\begin{equation} \label{eqn:so_linsys}
  \SoSys\colon \left\{
    \begin{aligned}
      (s^2 \BM + s\BD + \BK)\Bp(s) &= \Bg u(s), \\
      \Bz(s) & = \BC \Bp(s),
    \end{aligned}\right.
\end{equation}
where $\BM, \BD, \BK \in \C^{\nso \times \nso}$ are the mass, damping, and
stiffness matrices, while $\Bg \in \C^{\nso}$ and $\BC\in\C^{p\times \nso}$
are the linear input and output matrices.
The internal state $\Bp\colon \C \to \C^{\nso}$ of the plate system contains
$\nso=201\,900$ spatial coordinates that model the position of the plate in
the $x$, $y$, and $z$-directions. 
The RMS displacement of interest is recovered with $p = 27\,278$
entries of the state vector $\Bp$, considering only points on the inside of the
plate but not on the boundary.
This modeling approach is similar to that of~\eqref{eqn:linsys}.
For further details on the example, we refer the reader to~\cite{AumW23} and
the references therein.

Here, instead of the linear outputs in~\cref{eqn:so_linsys}, we consider 
directly the RMS error as the quantity of interest in a \emph{first-order}
linear quadratic output system of the form~\eqref{eqn:qosys}.
An equivalent system to that of~\eqref{eqn:so_linsys} is obtained by
rewriting the $\nso$ second-order equations in~\eqref{eqn:so_linsys} as
$n = 2\nso$ \emph{first-order} equations.
Therefore, we define the lifted state
\begin{equation*}
  \Bx = \begin{bmatrix} \Bp \\ s\Bp \end{bmatrix} \colon \C \to \Cn.
\end{equation*} 
Then, $\Bx$ is the state of a system~\eqref{eqn:qosys} where the system
matrices $\BE, \BA, \BQ \in \Cnn$ and $\Bb\in\Cn$ are defined in terms of the
system matrices in~\eqref{eqn:so_linsys} by
\begin{equation} \label{eqn:so_to_fo}
  \begin{aligned}
    \BE & = \begin{bmatrix} \BI_{\nso} & \Bzero \\ \Bzero & \BM \end{bmatrix}, &
    \BA & = \begin{bmatrix} \Bzero & \BI_{\nso} \\ -\BK & -\BD \end{bmatrix},\\
    \BQ & = \begin{bmatrix} \BC^{\trans}\BC & \Bzero \\ \Bzero & \Bzero
      \end{bmatrix}, &
    \Bb & =\begin{bmatrix} \Bzero \\ \Bg \end{bmatrix}.
\end{aligned}
\end{equation}
The resulting linear quadratic-output system of the form~\cref{eqn:qosys}
with the system matrices~\cref{eqn:so_to_fo} has the dimension $n = 403\,800$.
We consider this formulation of the system to construct suitable surrogates of
the form~\eqref{eqn:qosys_red} that can be used in place of the original
large-scale model for the frequency response analysis.

\section{Quadratic-output systems in time domain}%
\label{sec:lqo_review}

The description of linear quadratic-output systems in time domain as given in
the literature, e.g.,~\cite{VanM10, GosG22, BenGP21, ReiPGetal24}, strongly
resembles the frequency domain formulation~\cref{eqn:qosys} that we consider in
this work.
In the time domain, linear quadratic-output systems are typically represented as
\begin{equation} \label{eqn:qosys_td}
  \QOSystd\colon \left\{
    \begin{aligned}
      \BE\Bdxt(t) & = \BA\Bxt(t) +\Bb \ut(t), \quad \Bxt(0) = \Bzero_{n},\\
      \yt(t)^{2} & =\Bxt(t)^{\trans} \BQ \Bxt(t),
    \end{aligned}\right.
\end{equation}
with $\BA, \BE, \BQ \in \Rnn$, $\BQ$ symmetric, $\Bb \in \Rn$ and the
homogeneous initial condition $\Bxt(0)=\Bzero_{n}$.
Similar to~\cref{eqn:qosys}, the internal state of~\cref{eqn:qosys_td}
is denoted as $\Bxt\colon [0,\infty) \to \Rn$ and the external input as
$\ut\colon [0,\infty) \to \R$.
Despite bearing obvious similarities, the
systems~\cref{eqn:qosys,eqn:qosys_td} are not equivalent representations
connected via the Laplace transformation~\cite[Chap.~2.6]{DulP00}.
In fact, the nonlinear nature of the output equation in~\cref{eqn:qosys_td}
prevents the closed representation of the Laplace transform in the frequency
domain.
An alternative representation for the time domain systems~\cref{eqn:qosys_td}
in the Laplace domain is based on the idea of kernel representations and the
multivariate Laplace transform; see, e.g.,~\cite{DiaHGetal23}.
The resulting multivariate transfer function
$\Ht\colon \C \times \C \to \C$ of~\cref{eqn:qosys_td} is defined as
\begin{equation} \label{eqn:qotf_td}
  \Ht(s_{1}, s_{2}) = \Bb^{\trans} (s_{1} \BE - \BA)^{-\trans} \BQ
    (s_{2} \BE - \BA)^{-1} \Bb.
\end{equation}
In contrast, a transfer function of the frequency domain linear quadratic-output
system in~\cref{eqn:qosys} can be derived as follows:
At any point $s$ at which the matrix pencil $s\BE - \BA$ is nonsingular, the
state $\Bx(s)$ of the system in~\cref{eqn:qosys} is given explicitly by
\begin{equation*}
  \Bx(s) = (s\BE - \BA)^{-1} \Bb u(s).
\end{equation*}
Inserting $\Bx(s)$ into the output equation of~\cref{eqn:qosys} reveals the
input-to-output relationship of the system to be given as
\begin{align*}
  y(s) &= \left( (s\BE - \BA)^{-1} \Bb u(s)\right)^{\herm} \BQ
    (s\BE - \BA)^{-1} \Bb u(s)\\
  & = \Bb^{\herm} (s\BE - \BA)^{-\herm} \BQ (s\BE - \BA)^{-1} \Bb
    \lvert u(s) \rvert^{2}.
\end{align*}
The scalar-valued function $H\colon \C \to \C$, defined as
\begin{equation} \label{eqn:qotf}
  H(s) = \Bb^{\herm} (s\BE - \BA)^{-\herm} \BQ (s\BE - \BA)^{-1} \Bb,
\end{equation}
is the transfer function of the frequency domain linear qua\-dra\-tic-output
system~\cref{eqn:qosys}. 
Compared to $\Ht$ in~\cref{eqn:qotf_td}, there are some notable differences.
The transfer function~\cref{eqn:qotf} is univariate, and has a physical
interpretation in terms of the input-to-output mapping of~\cref{eqn:qosys} via
$y(s) = H(s) \lvert u(s) \rvert^{2}$.
Additionally, this transfer function~\cref{eqn:qotf} depends explicitly on
the complex conjugate of its argument, $\overline{s}$, and is thus
\emph{not analytic} on its domain, whereas the transfer
function~\cref{eqn:qotf_td} is analytic in both independent arguments.

The construction of reduced-order models for linear qua\-dra\-tic-output systems
has received increased consideration in recent years.
Proposed model reduction methods can broadly be categorized as approaches
that are based on energy functionals and balancing of the corresponding
system states~\cite{VanM10, PulA19, BenGP21, Pul23}, 
and the rational interpolation of transfer functions in the
frequency domain~\cite{VanVLetal12, GosA19, GosG22, DiaHGetal23, ReiPGetal23,
ReiPGetal24}.
Except for~\cite{VanM10, VanVLetal12}, all of these works consider solely
time domain quadratic-output systems~\cref{eqn:qosys_td}
and their corresponding multivariate transfer functions~\cref{eqn:qotf_td}.
The works~\cite{VanM10, VanVLetal12} do consider in parts also the frequency
domain system~\cref{eqn:qosys} but only via the intermediate representation
using a fully linear system of the form~\cref{eqn:linsys}, which as outlined
above has several disadvantages. 

Since we consider in this work the linear quadratic-output
system~\cref{eqn:qosys} in the frequency domain with the explicit transfer
function~\cref{eqn:qotf}, we primarily focus on interpolatory methods
that aim for an accurate approximation of this transfer function.
To this end, we will make use of the theory and ideas developed
in~\cite{GosA19, DiaHGetal23, ReiPGetal23} for time domain linear
quadratic-output systems as well as propose a new interpolatory framework for
the system class~\cref{eqn:qosys} that we consider in this work.


\section{Interpolation of linear quadratic-output systems}
\label{sec:interp_mor}

In the following, we present interpolation-based model order reduction for
frequency domain linear quadratic-output systems~\cref{eqn:qosys}.
First, we give an overview about projection-based model reduction
and interpolatory methods for the linear system case~\cref{eqn:linsys}.
Afterwards, we present new results that show how to enforce Lagrange and
Hermite interpolation conditions for~\cref{eqn:qotf} via projection of the
system matrices.


\subsection{Projection-based model order reduction}

Consider the linear quadratic-output system~\cref{eqn:qosys} and some given
reduction order $r \ll n$. 
In projection-based model order reduction, right and left basis matrices
$\BVr \in \Cnr$ and $\BWr \in \Cnr$ corresponding to $r$-dimensional
projection spaces $\mspan(\BVr)$ and $\mspan(\BWr)$ are chosen so that
the reduced-order model~\cref{eqn:qosys_red} is constructed as
\begin{equation} \label{eqn:proj}
  \begin{aligned}
    \BEr & = \BWr^{\herm} \BE \BVr, &
    \BAr & = \BWr{^\herm} \BA \BV ,\\
    \Bbr & = \BWr^{\herm} \Bb, &
    \BQr & = \BVr^{\herm} \BQ \BVr.
  \end{aligned}
\end{equation}
The corresponding transfer function of the resulting re\-duced-order
model~\cref{eqn:qosys_red} is given by
\begin{equation} \label{eqn:qotf_red}
  \Hr(s)= \Bbr^{\herm} (s\BEr - \BAr)^{-\herm} \BQr (s\BEr - \BAr)^{-1} \Bbr.
\end{equation}
Note that the same projection framework can be analogously applied to compute
reduced-order models of the time domain system~\cref{eqn:qosys_td}.
In this case, the transfer function of the reduced-order time domain system is
given by
\begin{equation} \label{eqn:qotf_td_red}
  \Hrt(s_{1}, s_{2}) = \Bbr^{\trans} (s_{1} \BEr - \BAr)^{-\trans} \BQr
    (s_{2} \BEr - \BAr)^{-1} \Bbr.
\end{equation}
In the purely linear setting~\cref{eqn:linsys}, the concept of
pro\-jec\-tion-based model order reduction is essentially the same.
The only difference is that in the output equation $\BCr = \BC \BVr$ replaces
$\BQr$ in~\cref{eqn:proj} such that the resulting reduced-or\-der transfer
function of the linear-output model becomes
\begin{equation} \label{eqn:lintf_red}
  \BGr(s)=\BCr(s \BEr - \BAr)^{-1}\Bbr.
\end{equation}
It is important to note that the reduced-order system is determined only
by the choice of the projection spaces $\mspan(\BVr)$ and $\mspan(\BWr)$
rather than the specific realization of the basis matrices $\BVr$ and
$\BWr$; cf.,~\cite[Sec.~3.3]{AntBG20}.
A classical choice to avoid ill-conditioned computations with the reduced-order
model is to replace any primitive constructions of $\BVr$ and $\BWr$ with
orthogonal matrices computed via, for example, a QR factorization or SVD.
In the special case that $\BA$ and $\BE$ are Hermitian or only one of the
two projection spaces can be computed, one typically sets $\BWr = \BVr$.
This setting is known in the literature as \emph{Galerkin projection}, while
the general case with $\BWr \neq \BVr$ and $\mspan(\BVr) \neq \mspan(\BWr)$ is
referred to as \emph{Petrov-Galerkin projection}.
In essence, any projection-based model order reduction method amounts to
choosing $\BVr$ and $\BWr$.
In the following sections, we consider the construction of the basis matrices
such that interpolation of the transfer function is enforced at selected points
in the complex plane.


\subsection{Interpolatory model reduction for linear systems}%
\label{ss:interp_linearmor}

In this section, we consider the interpolation-based model reduction of the
linear-output systems of the form~\eqref{eqn:linsys}.
In\-ter\-po\-la\-tion-based model reduction of linear systems is a well-studied
subject in the literature.
The theoretical foundation for projection-based interpolation has been laid
in~\cite{VilS87}, which was later developed into a numerically efficient
framework based on rational Krylov subspaces~\cite{Gri97}. 

Consider a linear system of the form~\cref{eqn:linsys} and given complex
interpolation points $\sigma_{1}, \ldots, \sigma_{k} \in \C$.
The idea of interpolatory methods is to construct a reduced-order linear model
with the corresponding transfer function $\BGr$ of the
form~\eqref{eqn:lintf_red} so that $\BGr$ interpolates the full-order
transfer function $\BG$ in~\cref{eqn:lintf} at the selected points, i.e.
\begin{equation} \label{eqn:lti_interp}
    \BG(\sigma_{i}) = \BGr(\sigma_{i}),
\end{equation}
for all $i = 1, \ldots, k$.
In the case of Hermite interpolation, additional conditions are set for the
derivatives of the transfer functions to satisfy
\begin{equation} \label{eqn:lti_hermite_interp}
  \frac{\operatorname{d}^{j_{i}}}{\operatorname{d}s^{j_{i}}} \BG(\sigma_{i})
  = \frac{\operatorname{d}^{j_{i}}}{\operatorname{d}s^{j_{i}}} \BGr(\sigma_{i}),
\end{equation}
for $i = 1, \ldots, k$ and some $1 \leq j_i\leq j$.
Because the full and reduced-order transfer functions $\BG$ and $\BGr$ assume
values in $\Cp$, enforcing interpolation of the entire transfer function can
quickly lead to reduced models of high orders when $p$ is large.
Recall that this is typically the case when rewriting the RMS error
measure~\cref{eqn:rms} into linear-output form.
A popular approach to remedy this is imposing interpolation conditions into
tangential directions:
Given the interpolation points $\sigma_{1}, \ldots, \sigma_{k} \in \C$ as well
as the direction vectors $\Bc_{1}, \ldots, \Bc_{k} \in \Cp$, the tangential
interpolation conditions associated with~\cref{eqn:lti_interp} are
\begin{equation*}
  \Bc_{i}^{\herm} \BG(\sigma_{i}) = \Bc_{i}^{\herm} \BGr(\sigma_{i}),
\end{equation*}
for $i = 1, \ldots, k$. 
Hermite tangential interpolation conditions analogous
to~\cref{eqn:lti_hermite_interp} can be enforced as well;
see~\cite[Thm.~3.3.2]{AntBG20} for more details.

The construction of appropriate interpolants can be done via projection.
Assume that the basis matrices $\BVr\in\Cnr$ and $\BWr\in\Cnr$ are chosen so
that
\begin{align*}
  (\sigma_{i} \BE - \BA)^{-1} \Bb & \in \mspan(\BVr) \quad\text{and}\\
  (\sigma_{i} \BE - \BA)^{-\herm} \BC^{\herm} \Bc & \in \mspan(\BWr)
\end{align*}
hold, for $i=1,\ldots, k$, then the reduced-order transfer function $\BGr$ that
is constructed via projection will be a tangential Hermite interpolant of
the full-order transfer function $\BG$ at the designated points and along
the tangential directions~\cite[Thm.~3.3.1]{AntBG20}.

The choice of interpolation points is a crucial element in 
interpolatory model-order reduction, which ultimately determines the
performance of the reduced-order model.
To select good or even optimal interpolation points, different approximation
error measures have been used in the literature.
The most common one is the $\CH_{2}$-norm, for which locally optimal
reduced-order models can be constructed via tangential Hermite
interpolation~\cite{MeiL67, GugAB08}.
The associated interpolation points are the mirror images of the eigenvalues of
the reduced-order matrix pencil $\lambda\BEr - \BAr$ and the tangential
directions are the corresponding transfer function residues.
The iterative rational Krylov algorithm (IRKA) is an efficient approach to
compute these optimal interpolation points and tangential 
directions via iterative updates~\cite{GugAB08}.
A different error measure is the $\CL_{\infty}$-norm, for which typically
greedy procedures are used that select interpolation points iteratively as
maxima of the error transfer function $\BG - \BGr$; see, for
example,~\cite{AliBMetal20, AumW23}.


\subsection{Interpolation in the quadratic-output setting}

The aforementioned interpolation of linear systems can be extended to the time
domain quadratic-output systems~\cref{eqn:qosys_td}.
For the transfer function~\cref{eqn:qotf_td}, given a set of interpolation
points $\sigma_{1}, \ldots, \sigma_{k} \in \C$, the multivariate interpolation
problem is to find $\Hrt$ such that
\begin{equation*}
  \Ht(\sigma_{i}, \sigma_{j}) = \Hrt(\sigma_{i}, \sigma_{j})
\end{equation*}
holds, for $i, j = 1, \ldots, k$. 
As in the linear setting, these interpolation conditions can be enforced via
projection using rational Kyrlov subspaces; see~\cite[Cor.~1]{DiaHGetal23}.
The work~\cite{ReiPGetal23} shows that $\CH_{2}$-optimal interpolants for
time domain systems of the form~\cref{eqn:qosys_td} satisfy the following
Hermite interpolation conditions
\begin{equation*}
  \Ht(-\overline{\lambda}_{i}, -\overline{\lambda}_{j})  =
    \Hrt(-\overline{\lambda}_{i}, -\overline{\lambda}_{j}),
\end{equation*}
and
\begin{align*}
  & \sum_{k = 1}^{r} \left(\qr_{i,k} \frac{\partial}{\partial s_{1}}
    \Ht(-\overline{\lambda}_{i}, -\overline{\lambda}_{k}) +
    \qr_{k,i}\frac{\partial}{\partial s_{2}}
    \Ht(-\overline{\lambda}_{k}, -\overline{\lambda}_{i})\right) \\ 
  &\quad = \sum_{k=1}^{r} \left(\qr_{i,k}
    \frac{\partial}{\partial s_{1}}
    \Hrt(-\overline{\lambda}_{i}, -\overline{\lambda}_{k}) + \qr_{k,i}
    \frac{\partial}{\partial s_{2}}
    \Hrt(-\overline{\lambda}_{k}, -\overline{\lambda}_{i})\right),
\end{align*}
for $i, j = 1, \ldots, r$, where $\lambda_{1}, \ldots, \lambda_{r}\in\C$ are
the eigenvalues of the reduced-order matrix pencil $\lambda\BEr - \BAr$
in~\cref{eqn:proj}, and $\qr_{i,j} \in \C$ denotes the $(i, j)$-th
entry of $\BQr \in \Crr$ in an appropriate basis.
Based on these interpolation conditions, a generalization of the IRKA approach
for systems of the form~\cref{eqn:qosys_td} is proposed in~\cite{ReiPGetal23}.
Due to the similarities between~\cref{eqn:qosys,eqn:qosys_td}
as well as their transfer functions~\cref{eqn:qotf_td,eqn:qotf},
we expect this IRKA-like approach to be well suited for constructing accurate
projection spaces for the frequency domain quadratic-output
system~\cref{eqn:qosys}, too.

In the following, we develop a new projection-based interpolation theory for
the frequency domain quadratic-output system~\cref{eqn:qosys}.
As previously mentioned, the transfer function~\cref{eqn:qotf} is not
analytic such that Hermite interpolation requires special care.
However, in most instances, the interpolation points for model order reduction
are selected only from the imaginary axis~\cite{AliBMetal20, AumW23}.
Considering $z \in \i\R$, the transfer function~\cref{eqn:qotf} simplifies to
\begin{equation} \label{eqn:qotfimag}
  H(z) = \Bb^{\herm} (-z \BE^{\herm} - \BA^{\herm})^{-1} \BQ
    (z \BE - \BA)^{-1} \Bb.
\end{equation}
In  this case, the transfer function~\cref{eqn:qotfimag} is differentiable
along the imaginary axis, which leads to the following theorem.

\begin{theorem} \label{thm:interp}
  Let $H\colon \i\R \to \C$ be the transfer function~\cref{eqn:qotf} and
  $\Hr\colon \i\R \to \C$ be the reduced-order transfer
  function~\cref{eqn:qotf_red} obtained via projection~\cref{eqn:qosys_red},
  both restricted to the imaginary axis.
  Consider interpolation points $z_{1}, \ldots, z_{k} \in\i\R$, in which
  $H$ and $\Hr$ can be evaluated, and the basis matrices $\BVr \in \Cnr$ and
  $\BWr \in \Cnr$ satisfy
  \begin{align} \label{eqn:lqo_Vr}
    (z_{i} \BE - \BA)^{-1} \Bb & \in \mspan(\BVr),\\
    \label{eqn:lqo_Wr}
    (z_{i} \BE - \BA)^{-\herm} \BQ (z_{i} \BE - \BA)^{-1} \Bb &
      \in \mspan (\BWr), 
  \end{align}
  for $i = 1, \ldots, k$.
  Then, for $i = 1, \ldots, k$, the following interpolation conditions hold:
  \begin{equation*}
    H(z_{i}) = \Hr(z_{i}) \quad\text{and}\quad
    \frac{\operatorname{d}}{\operatorname{d} z}H(z_{i})
      = \frac{\operatorname{d}}{\operatorname{d} z}\Hr(z_{i}).
  \end{equation*}
\end{theorem}
\medskip
\begin{proof}
  Set $\CK(z) = z \BE - \BA$,
  $\CKr(s) = s\BEr-\BAr$, $\Bv_{i} = \CKr(z_{i})^{-1} \Bbr$
  and $\Bw_{i} = \CKr(z_{i})^{-\herm} \BQr \CKr(z_{i})^{-1} \Bbr$.
  By construction of $\BV$ and $\BW$, the following two projection identities
  hold
  \begin{align*}
    \BV \Bv_{i} & = \CK(z_{i})^{-1} \Bb, \quad \BW \Bw_{i} = \CK(z_{i})^{-\herm}
      \BQ \CK(z_{i})^{-1} \Bb;
  \end{align*}
  see, e.g.,~\cite{Wer21} for details on the projectors.
  It follows that
  \begin{equation*}
    \Hr(z_{i}) = \Bv_{i}^{\herm} \BQr \Bv_{i}
      = \Bv_{i}^{\herm} \BV^{\herm} \BQ  \BV \Bv_{i} = H(z_{i})
  \end{equation*}
  holds, for $i = 1, \ldots, k$.
  The derivative of the reduced-order transfer function is given via
  \begin{equation*}
    \frac{\operatorname{d}}{\operatorname{d} z} \Hr(z_{i}) =
      \Bv_{i}^{\herm} \BEr^{\herm} \Bw_{i} - \Bw_{i}^{\herm} \BEr \Bv_{i}.
  \end{equation*}
  Applying both the identities above yields
  \begin{equation*}
    \frac{\operatorname{d}}{\operatorname{d} z} \Hr(z_{i})
    = \Bv_{i}^{\herm} \BV^{\herm} \BE^{\herm} \BW \Bw_{i} -
      \Bw_{i}^{\herm} \BW^{\herm} \BE \BV \Bv_{i}
    = \frac{\operatorname{d}}{\operatorname{d} z} H(z_{i}),
  \end{equation*}
  for $i = 1, \ldots, k$, which concludes the proof.
\end{proof}

The Lagrange interpolation condition in \Cref{thm:interp} is similar
to the result in~\cite[Cor.~1]{DiaHGetal23} for multivariate transfer
functions~\cref{eqn:qotf_td}.
However, \Cref{thm:interp} is specifically tailored for~\cref{eqn:qotf},
which we consider in this paper, and provides new Hermite interpolation results.
Note that due to the structure of the primitive bases vectors
in~\cref{eqn:lqo_Vr,eqn:lqo_Wr}, the basis matrices $\BVr \in \Cnr$
and $\BWr \in \Cnr$ can be computed efficiently by solving linear systems
of equations.
If the Hermite interpolation conditions are not required,
$\BWr$ can also be used to match additional Lagrange interpolation conditions by
mimicking the structure of $\BVr$ at additional points.

\section{Numerical experiments}%
\label{sec:numerics}

In the following, we apply the interpolatory model-order reduction
schemes from \Cref{sec:interp_mor} to the vibrating plate model in
\Cref{sec:plate_tva}.
The experiments reported here have been executed on the \texttt{Gauss} compute
server of the Department of Mathematics at Virginia Tech.
This machine is equipped with 4 Intel(R) Xeon(R) CPU E7-4890 v2 processors
running at 3.17\,GHz and 1\,TB total main memory.
It is running on Ubuntu 22.04.3 LTS and uses MATLAB 23.2.0.2365128 (R2023b).
The source code for the numerical experiments and the computed results are
available at~\cite{supReiW24a}.


\subsection{Experimental setup}%
\label{ss:exp_setup}

For the model reduction of frequency domain quadratic-output
systems~\cref{eqn:qosys}, we propose and compare five different
interpolation-based approaches as described in the following:
\begin{description}
  \item[\irka{}] is the time domain quadratic-output IRKA that we
    apply to the matrices of the frequency domain system as described in
    \Cref{sec:interp_mor};
  \item[\intinfv{}] performs interpolation following
    \Cref{thm:interp} using
    a Galerkin projection based on~\cref{eqn:lqo_Vr} and selecting
    interpolation points greedily from a pre-sampled basis;
  \item[\intinfvw{}] is the same as \intinfv{} but with a Petrov-Galerkin
    projection using~\cref{eqn:lqo_Vr,eqn:lqo_Wr};
  \item[\intavgv{}] uses a Galerkin projection based on a truncated
    approximation of a pre-sampled interpolation basis as in
    \Cref{thm:interp} and~\cref{eqn:lqo_Vr} computed by the pivoted QR
    decomposition,
  \item[\intavgvw{}] is the same as \intavgv{} but with a Petrov-Galerkin
    projection using two pre-sampled bases
    following~\cref{eqn:lqo_Vr,eqn:lqo_Wr}.
\end{description}
The latter four of these methods ($\mathsf{Int_{*}}$) rely on a pre-sampled
basis, for which we compute and save the solutions
to~\cref{eqn:lqo_Vr,eqn:lqo_Wr} in $250$ linearly equidistant points in
the interval $\i [1, 2\pi \cdot 251] \subset \i\R$.

To visibly compare different approximations, we use the pointwise relative
approximation errors of the transfer functions as
\begin{equation} \label{eqn:relerr}
  \relerr(\omega) = \frac{\lvert H(2\pi\i \cdot \omega) -
    \Hr(2\pi\i \cdot \omega) \rvert}{\lvert H(2\pi\i \cdot \omega)\rvert},
\end{equation}
for the frequencies $\omega \in [0, 250]$\,Hz, in plots alongside the magnitude
of the transfer functions.
Additionally, we score the performance of the proposed methods based on two
relative system error measures.
To compare the average performance of the computed reduced models over
the frequency range $\i [0, 2\pi \cdot 250] \subset \i\R$, we compute an
approximation to the relative $\CH_{2}$ error as
\begin{equation} \label{eqn:H2relerr}
  \relerr_{\CH_{2}} = \frac{\sum_{i = 0}^{500} \lvert
    H(2 \pi \i \cdot \omega_{i}) - \Hr(2 \pi \i \cdot \omega_{i}) \rvert}%
    {\sum_{i = 0}^{500} \lvert H(2\pi\i \cdot \omega_{i})\rvert},
\end{equation}
where $\Omega = \{\omega_{i}\}_{i = 1}^{500}$ is the collection of $500$
linearly equidistant points in the interval
$[0, 250] \subset \R$.
Also, to compare the worst case performance of the computed reduced models over
the frequency range $[0, 250]$\,Hz, we compute the following approximation to
the relative $\CH_{\infty}$ error via
\begin{equation} \label{eqn:Hinftyrelerr}
  \relerr_{\CH_{\infty}} =\frac{\max_{\omega_{i} \in \Omega} \lvert
    H(2 \pi \i \cdot \omega_{i}) - \Hr(2 \pi \i \cdot \omega_{i}) \rvert}%
    {\max_{\omega_{i} \in \Omega} \lvert H(2 \pi \i \cdot \omega_{i})\rvert},
\end{equation}
where $\Omega$ is the set of $500$ frequency points as defined above.


\subsection{Discussion}%
\label{sec:discussion}

\begin{figure*}[t!]
  \centering
  \vspace{.5\baselineskip}
  \begin{subfigure}[b]{.49\linewidth}
    \raggedleft
  \tikzexternalenable%
  \tikzsetnextfilename{r50_mag}%
  \input{graphics/r50_mag.tikz}%
  \tikzexternaldisable%
\\
  \tikzexternalenable%
  \tikzsetnextfilename{r50_err}%
  \input{graphics/r50_err.tikz}%
  \tikzexternaldisable%

    \caption{Transfer function magnitudes and relative errors for $r = 50$.}
    \label{fig:50}
  \end{subfigure}%
  \hfill%
  \begin{subfigure}[b]{.49\linewidth}
    \raggedleft
  \tikzexternalenable%
  \tikzsetnextfilename{r75_mag}%
  \input{graphics/r75_mag.tikz}%
  \tikzexternaldisable%
\\[-.4\baselineskip]
  \tikzexternalenable%
  \tikzsetnextfilename{r75_err}%
  \input{graphics/r75_err.tikz}%
  \tikzexternaldisable%

    \caption{Transfer function magnitudes and relative errors for $r = 75$.}
    \label{fig:r75}
  \end{subfigure}

  \vspace{.5\baselineskip}
  \tikzexternalenable%
  \tikzsetnextfilename{legend}%
  \input{graphics/legend.tikz}%
  \tikzexternaldisable%

  \caption{Frequency response results for reduced-order models of orders
    $r = 50$ and $r=75$:
    All proposed methods provide reasonably accurate approximations of the
    original system behavior.
    The interpolation methods that take the output equation into account
    provide visibly more accurate results than the classical Galerkin-based
    methods.
    \irka{} performs multiple orders of magnitude worse for $r=75$ compared to
    the other methods due to the method not converging in the prescribed amount 
    of iteration steps.}
  \label{fig:numerics}
\end{figure*}
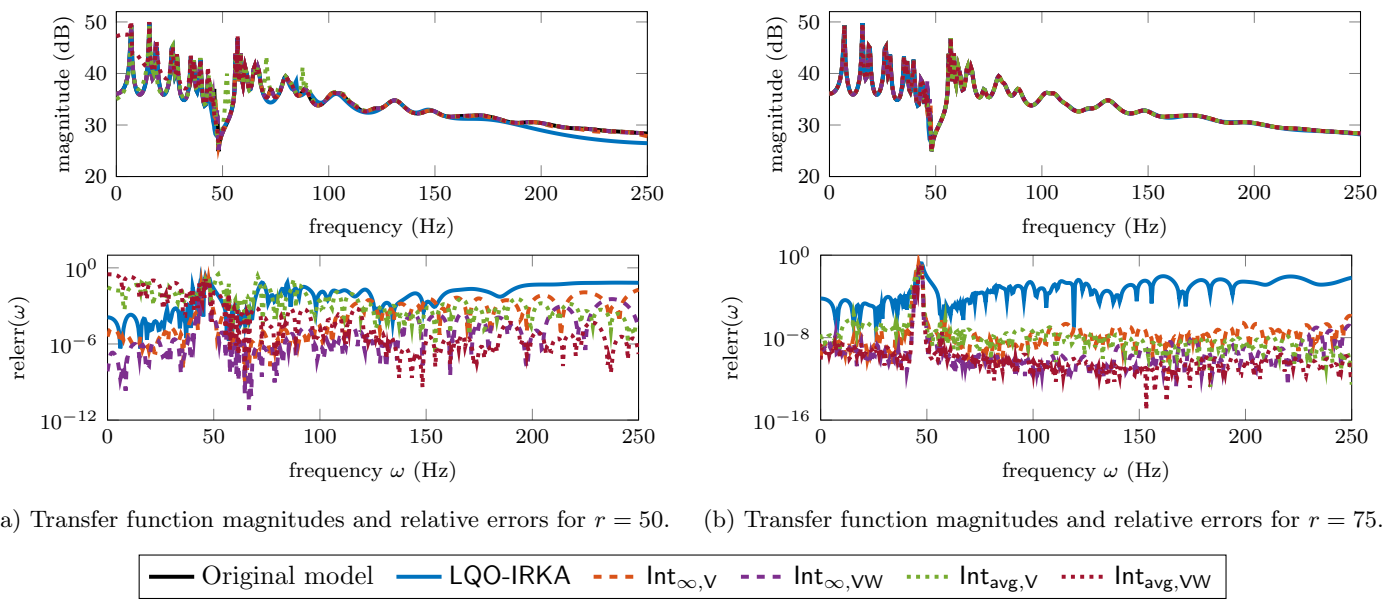

We computed reduced-order models of orders $r = 25, 50, 75$ and $100$ using the
five approaches described in the previous section.
The performance of the computed reduced models for orders $r = 50$ and $r = 75$
in the considered frequency range is shown \Cref{fig:numerics} in terms of the
transfer function magnitudes and the pointwise relative approximation
errors~\cref{eqn:relerr}.
The relative approximate system norm errors for each method and all four
reduction orders $r$ are recorded in~\Cref{tab:relH2,tab:relHinfty}.

\begin{table}[t]
  \centering
  \caption{Relative $\CH_{2}$ errors according
    to~\eqref{eqn:H2relerr} for reduced-order models of size
    $r = 25, 50, 75, 100$.
    The smallest error for each order is highlighted in \textbf{boldface}.}
  \label{tab:relH2}
  \vspace{.5\baselineskip}
  
  \resizebox{\linewidth}{!}{%
  \begin{tabular}{lrrrr}
    \hline\noalign{\smallskip}
      & \multicolumn{1}{c}{$r = 25$}
      & \multicolumn{1}{c}{$r = 50$}
      & \multicolumn{1}{c}{$r = 75$}
      & \multicolumn{1}{c}{$r = 100$}  \\
    \noalign{\smallskip}\hline\noalign{\smallskip}
    \irka
      & $1.071\texttt{e-}1$
      & $1.973\texttt{e-}2$
      & $2.529\texttt{e-}3$
      & $1.017\texttt{e-}3$ \\
    \intinfv
      & $1.119\texttt{e-}1$
      & $2.466\texttt{e-}3$
      & $7.357\texttt{e-}4$
      & $9.801\texttt{e-}4$ \\
    \intinfvw
      & $\boldsymbol{7.013\texttt{e-}2}$
      & $\boldsymbol{1.005\texttt{e-}3}$
      & $\boldsymbol{4.616\texttt{e-}4}$
      & $\boldsymbol{3.841\texttt{e-}4}$ \\
    \intavgv
      & $1.020\texttt{e-}1$
      & $1.145\texttt{e-}2$
      & $5.722\texttt{e-}4$
      & $5.015\texttt{e-}4$ \\
    \intavgvw
      & $2.364\texttt{e-}1$
      & $1.783\texttt{e-}1$
      & $4.674\texttt{e-}4$
      & $5.519\texttt{e-}4$ \\
    \noalign{\smallskip}\hline\noalign{\smallskip}
  \end{tabular}}
\end{table}

\begin{table}[t]
  \centering
  \caption{Relative $\CH_{\infty}$ errors according
    to~\eqref{eqn:Hinftyrelerr} for reduced-order models of size
    $r = 25, 50, 75, 100$.
    The smallest error for each order is highlighted in \textbf{boldface}.}
  \label{tab:relHinfty}
  \vspace{.5\baselineskip}
  
  \resizebox{\linewidth}{!}{%
  \begin{tabular}{lrrrr}
    \hline\noalign{\smallskip}
      & \multicolumn{1}{c}{$r = 25$}
      & \multicolumn{1}{c}{$r = 50$}
      & \multicolumn{1}{c}{$r = 75$}
      & \multicolumn{1}{c}{$r = 100$}  \\
    \noalign{\smallskip}\hline\noalign{\smallskip}
    \irka
      & $2.169\texttt{e-}1$
      & $1.539\texttt{e-}1$
      & $1.066\texttt{e-}1$
      & $6.436\texttt{e-}2$ \\
    \intinfv
      & $1.769\texttt{e-}1$
      & $\boldsymbol{8.872\texttt{e-}2}$
      & $9.475\texttt{e-}2$
      & $1.168\texttt{e-}1$ \\
    \intinfvw
      & $\boldsymbol{1.569\texttt{e-}1}$
      & $9.291\texttt{e-}2$
      & $\boldsymbol{5.030\texttt{e-}2}$
      & $\boldsymbol{6.076\texttt{e-}2}$ \\
    \intavgv
      & $5.159\texttt{e-}1$
      & $1.918\texttt{e-}1$
      & $6.368\texttt{e-}2$
      & $9.623\texttt{e-}2$ \\
    \intavgvw
      & $7.023\texttt{e-}1$
      & $2.254\texttt{e-}1$
      & $1.058\texttt{e-}1$
      & $8.828\texttt{e-}2$ \\
    \noalign{\smallskip}\hline\noalign{\smallskip}
    \end{tabular}}
\end{table}

We observe that all compared methods provide reasonably accurate approximations
for all considered orders of reduction $r = 25, 50, 75, 100$.
For the $r = 50$ reduced-order model in~\Cref{fig:numerics}, it is evident that
the errors of \intinfv{} and \intinfvw{} are multiple orders of magnitude
smaller for lower frequencies compared to the other three methods.
Similar behavior could be observed for the $r = 25$ reduced model; see the
accompanying code package for these results~\cite{supReiW24a} as well as for
$r = 100$.
This comes from the use of the absolute $\CH_{\infty}$ error measure in the
greedy approaches.
On the other hand, \intavgv{} and \intavgvw{} provide a better approximation
for higher frequencies through the basis approximation since all frequencies
are considered to be about equally important here.
This leads to less accurate approximations for lower frequencies, where a
lot of the dominant transfer function behavior happens.
\irka{} exhibits an overall similar approximation behavior.
Note that \intinfvw{} and \intavgvw{}, which include information about the
nonlinear output equation, perform visibly better than their Galerkin
counterparts that only account for the linear input-to-state equation.
This can also be seen in~\Cref{tab:relH2,tab:relHinfty}.
The relative error values of \intinfvw{} and \intavgvw{} are smaller than their
Galerkin counterparts for almost every order of reduction.
For $r=75$ and $r = 100$, the \irka{} iteration did not converge in the
prescribed number of steps such that its approximation accuracy cannot keep
up with the other approaches, which appear to have reached the smallest
possible error for most of the frequency interval.
We note that for any order, the frequency $48$\,Hz is difficult to match
due to the use of the TVAs; see~\cite{AumW23}.
This could be resolved by additionally enforcing interpolation at this
particular point.
Lastly, for every order of reduction and in each metric, the greedy
sampling-based methods \intinfv{} and \intinfvw{} exhibit the smallest
magnitude errors, with \intinfvw{} performing better in almost every instance. 
This holds true even in the relative $\CH_{2}$ error
metric~\eqref{eqn:H2relerr}, despite the expectation that the averaging-based
approaches may perform better here.

\section{Conclusions}%
\label{sec:conclusions}

In this work, we presented model order reduction approaches for systems with
root mean squared error measures in the frequency domain.
We provided a new theoretical framework for Hermite interpolation of such
systems and showed in numerical experiments that interpolation methods based on
this new theory outperform classical approaches that rely on the linearization
of the measured quantities.
The \irka{} method transcribed from time domain quadratic-output systems
performed comparably well when it converged, which leads to the question if one
can design a similar optimization procedure for the frequency domain systems
that we considered in this work.


\section*{Acknowledgments}%
\addcontentsline{toc}{section}{Acknowledgments}

The authors would like to thank Serkan Gugercin for inspiring discussions.


\addcontentsline{toc}{section}{References}
\bibliographystyle{plainurl}
\bibliography{bibtex/myref}
  
\end{document}

%% file: graphics/plate.tikz
\begingroup%
  \makeatletter%
  \providecommand\color[2][]{%
    \errmessage{(Inkscape) Color is used for the text in Inkscape, but the package 'color.sty' is not loaded}%
    \renewcommand\color[2][]{}%
  }%
  \providecommand\transparent[1]{%
    \errmessage{(Inkscape) Transparency is used (non-zero) for the text in Inkscape, but the package 'transparent.sty' is not loaded}%
    \renewcommand\transparent[1]{}%
  }%
  \providecommand\rotatebox[2]{#2}%
  \newcommand*\fsize{\dimexpr\f@size pt\relax}%
  \newcommand*\lineheight[1]{\fontsize{\fsize}{#1\fsize}\selectfont}%
  \ifx\svgwidth\undefined%
    \setlength{\unitlength}{403.2683189bp}%
    \ifx\svgscale\undefined%
      \relax%
    \else%
      \setlength{\unitlength}{\unitlength * \real{\svgscale}}%
    \fi%
  \else%
    \setlength{\unitlength}{\svgwidth}%
  \fi%
  \global\let\svgwidth\undefined%
  \global\let\svgscale\undefined%
  \makeatother%
  \begin{picture}(1,0.3474168)%
    \lineheight{1}%
    \setlength\tabcolsep{0pt}%
    \put(0,0){\includegraphics[width=\unitlength,page=1]{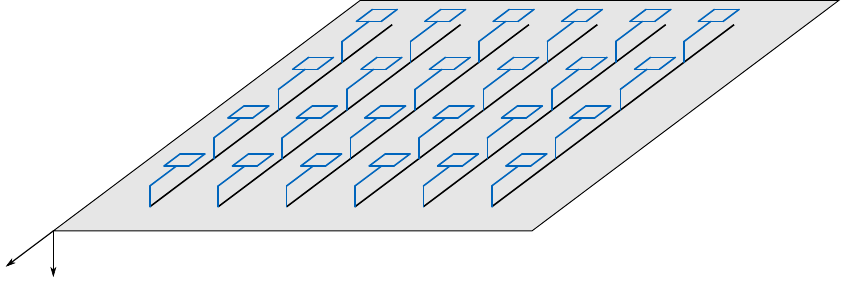}}%
    \put(0.07543807,0.02816779){\color[rgb]{0,0,0}\makebox(0,0)[lt]{\begin{minipage}{0.05313725\unitlength}\raggedright \Large $z$\end{minipage}}}%
    \put(0.69591664,0.10486813){\color[rgb]{0,0,0}\makebox(0,0)[lt]{\begin{minipage}{0.05313725\unitlength}\raggedright \Large $x$\end{minipage}}}%
    \put(-0.00017981,0.01855207){\color[rgb]{0,0,0}\makebox(0,0)[lt]{\begin{minipage}{0.05313725\unitlength}\raggedright \Large $y$\end{minipage}}}%
    \put(0,0){\includegraphics[width=\unitlength,page=2]{graphics/plate_sketch_svg-tex.pdf}}%
  \end{picture}%
\endgroup%

%% file: graphics/r50_mag.tikz
\begin{tikzpicture}[font = \plotfontsize]
  \pgfplotstableread{graphics/data/r50_mag.dat}\tableINPUT
  
  \begin{axis}[%
    width  = .775\linewidth,
    height = .085\textheight,
    scale only axis,
    xmin = 0,
    xmax = 250,
    ymin = 20,
    ymax = 52,
    xminorticks = false,
    yminorticks = false,
    xlabel = {frequency (Hz)},
    ylabel = {magnitude (dB)},
    ylabel style   = {yshift = -.3em},
    xlabel style   = {yshift = .3em},
    scaled x ticks = false,
    x tick label style = {/pgf/number format/1000 sep={\,}},
    y tick label style = {/pgf/number format/1000 sep={\,}},
    cycle list name    = plotlist
  ]
  
    \foreach \y in {1, 2, ..., 6}{
      \addplot+ table[x index = 0, y index = \y] {\tableINPUT};
    }
  \end{axis}
\end{tikzpicture}

%% file: graphics/r50_err.tikz
\begin{tikzpicture}[font = \plotfontsize]
  \pgfplotstableread{graphics/data/r50_err.dat}\tableINPUT
  
  \begin{semilogyaxis}[%
    width  = .775\linewidth,
    height = .085\textheight,
    scale only axis,
    xmin = 0,
    xmax = 250,
    ymin = 1.0e-12,
    ymax = 1.0e+1,
    xminorticks = false,
    yminorticks = false,
    xlabel = {frequency $\omega$ (Hz)},
    ylabel = {$\relerr(\omega)$},
    ylabel style   = {yshift = -.3em},
    xlabel style   = {yshift = .3em},
    scaled x ticks = false,
    x tick label style = {/pgf/number format/1000 sep={\,}},
    y tick label style = {/pgf/number format/1000 sep={\,}},
    cycle list name    = plotlist
  ]

  \pgfplotsset{cycle list shift = 1}
  
    \foreach \y in {1, 2, ..., 5}{
      \addplot+ table[x index = 0, y index = \y] {\tableINPUT};
    }
  \end{semilogyaxis}
\end{tikzpicture}

%% file: graphics/r75_mag.tikz
\begin{tikzpicture}[font = \plotfontsize]
  \pgfplotstableread{graphics/data/r75_mag.dat}\tableINPUT
  
  \begin{axis}[%
    width  = .775\linewidth,
    height = .085\textheight,
    scale only axis,
    xmin = 0,
    xmax = 250,
    ymin = 20,
    ymax = 52,
    xminorticks = false,
    yminorticks = false,
    xlabel = {frequency (Hz)},
    ylabel = {magnitude (dB)},
    ylabel style   = {yshift = -.3em},
    xlabel style   = {yshift = .3em},
    scaled x ticks = false,
    x tick label style = {/pgf/number format/1000 sep={\,}},
    y tick label style = {/pgf/number format/1000 sep={\,}},
    cycle list name    = plotlist
  ]
  
    \foreach \y in {1, 2, ..., 6}{
      \addplot+ table[x index = 0, y index = \y] {\tableINPUT};
    }
  \end{axis}
\end{tikzpicture}

%% file: graphics/r75_err.tikz
\begin{tikzpicture}[font = \plotfontsize]
  \pgfplotstableread{graphics/data/r75_err.dat}\tableINPUT
  
  \begin{semilogyaxis}[%
    width  = .775\linewidth,
    height = .085\textheight,
    scale only axis,
    xmin = 0,
    xmax = 250,
    ymin = 1.0e-16,
    ymax = 1.0e+0,
    xminorticks = false,
    yminorticks = false,
    xlabel = {frequency $\omega$ (Hz)},
    ylabel = {$\relerr(\omega)$},
    ylabel style   = {yshift = -.3em},
    xlabel style   = {yshift = .3em},
    scaled x ticks = false,
    x tick label style = {/pgf/number format/1000 sep={\,}},
    y tick label style = {/pgf/number format/1000 sep={\,}},
    cycle list name    = plotlist
  ]

  \pgfplotsset{cycle list shift = 1}
  
    \foreach \y in {1, 2, ..., 5}{
      \addplot+ table[x index = 0, y index = \y] {\tableINPUT};
    }
  \end{semilogyaxis}
\end{tikzpicture}

%% file: graphics/legend.tikz
\begin{tikzpicture}
  \begin{axis}[%
    hide axis,
    width  = 1mm,
    height = 1mm,
    scale only axis,
    xmin = 0,
    xmax = 1,
    ymin = 0,
    ymax = 1,
    legend columns = 6, 
    legend style   = {
      at     = {(0,0)},
      anchor = center,
      /tikz/every even column/.append style = {column sep = 0.2cm}},
    legend cell align  = {left},
    clip mode          = individual,
    cycle list name    = plotlist]

    \foreach \y in {1, 2, ..., 6}{
      \addplot+ coordinates{ (0, 0) };
    }
    
    \addlegendentry{Original model}
    \addlegendentry{\irka{}}
    \addlegendentry{\intinfv{}}
    \addlegendentry{\intinfvw{}}
    \addlegendentry{\intavgv{}}
    \addlegendentry{\intavgvw{}}
  \end{axis}
\end{tikzpicture}